\documentclass[a4paper,10pt]{amsart}
\usepackage{amssymb}

\numberwithin{equation}{section}
\newtheorem{theorem}{Theorem}[section]

\newtheorem{lem}[theorem]{Lemma}

\newtheorem{rem}[theorem]{Remark}
\theoremstyle{remark}

\newcommand{\R}{\mathbb{R}}

\author[H.~Orf]{Hajer Orf}
\address{Laboratory of Mathematics and Applications, College of Sciencs, Gabes University. Tunisia
}
\email{\sl hajerorf17@gmail.com}

\title[Long time decay  for global solutions... ]
{Long time decay  for global solutions to the Navier-Stokes equations in Sobolev--Gevrey spaces}
 
\begin{document}
\begin{abstract}
We study the three-dimensional incompressible Navier--Stokes equations in the critical Sobolev--Gevrey space $\dot{H}^{1/2}_{a,1}(\mathbb{R}^3)$, $a>0$. We first prove local well-posedness, together with a global existence criterion and a blow-up criterion for the maximal solution. Our main result is that every global solution $u\in C\big([0,+\infty),\dot{H}^{1/2}_{a,1}(\mathbb{R}^3)\big)$ satisfies $\|u(t)\|_{\dot{H}^{1/2}_{a,1}}\to0$ as $t\to+\infty$. The proof combines the classical decay of global solutions in $\dot{H}^{1/2}$ with an analyticity estimate for small data in $\dot{H}^{1/2}$ and Fourier-analytic product laws in Sobolev--Gevrey spaces. As consequences, we obtain the decay rate $\|u(t)\|_{\dot{H}^s}=o\big(t^{-\frac{s-1/2}{2}}\big)$ for every $s\ge\frac12$, and the stability of global solutions under small perturbations of the initial data in $\dot{H}^{1/2}_{a,1}$.
\end{abstract}
 

\maketitle
\tableofcontents

 
\section{Introduction}
The study of the long-time behavior of solutions to the three-dimensional incompressible Navier--Stokes equations remains one of the central problems in mathematical fluid dynamics. Despite significant progress over the past decades, understanding how viscous incompressible flows evolve and dissipate energy over time continues to pose deep analytical challenges. In particular, describing the decay properties of global solutions in critical or near-critical functional spaces provides essential insight into the mechanisms of stability and regularity for the Navier–Stokes dynamics. Establishing such decay results not only refines our understanding of the asymptotic regime of viscous flows but also contributes to the broader program of characterizing the global behavior of solutions in physically relevant spaces.
We consider the three-dimensional incompressible Navier--Stokes equations
\begin{equation}\label{NS}
\begin{cases}
\partial_t u-\Delta u+u\cdot\nabla u=-\nabla p, & (t,x)\in\mathbb{R}^+\times\mathbb{R}^3,\\
\operatorname{div}u=0, & (t,x)\in\mathbb{R}^+\times\mathbb{R}^3,\\
u(0,x)=u^0(x), & x\in\mathbb{R}^3,
\end{cases}
\tag{NS}
\end{equation}
where $u=(u_1,u_2,u_3)$ is the velocity field of the fluid, $p$ is the
pressure, and $u^0=(u^0_1,u^0_2,u^0_3)$ is the initial velocity. The
viscosity coefficient plays no role in our arguments, and we set it equal to
one. Since $\operatorname{div}u=0$, the pressure is determined by the velocity
through $p=(-\Delta)^{-1}\partial_i\partial_j(u_iu_j)$.
 
Let $\mathbb{P}$ be the Leray projector onto divergence-free vector fields,
whose Fourier symbol is $\delta_{jk}-\xi_j\xi_k/|\xi|^2$. Applying $\mathbb{P}$
to \eqref{NS} and using Duhamel's formula, we write the system in the
equivalent integral form
\begin{equation}\label{NS-integral}
u(t)=e^{t\Delta}u^0+B(u,u)(t),\qquad
B(u,v)(t):=-\int_0^t e^{(t-s)\Delta}\,\mathbb{P}\operatorname{div}(u\otimes v)(s)\,ds .
\end{equation}
 
The well-posedness problem for the Navier--Stokes equations has been extensively studied in various critical functional settings. Fujita and Kato \cite{FK} proved local existence for arbitrary data, and global existence for small data, in $\dot{H}^{1/2}(\mathbb{R}^3)$. Lei and Lin \cite{ZF} obtained global mild solutions for small data in a space defined through the Fourier transform, and Bae \cite{BAE} proved that these solutions are analytic in time; the analyticity estimate of Lemma~\ref{theo3} below is inspired by \cite{BAE}. Sobolev--Gevrey spaces $\dot{H}^{s}_{a,\sigma}$, which measure analytic ($\sigma=1$) or Gevrey ($\sigma>1$) regularity with the weight $e^{a|\xi|^{1/\sigma}}$, were used in \cite{BJ1} to study blow-up criteria.\\
 
Concerning the long-time behavior of global solutions, Gallagher, Iftimie and Planchon \cite{GI} proved that if $u\in C(\mathbb{R}^+,\dot{H}^{1/2})$ is a global solution of the three-dimensional Navier--Stokes equations, then
$$\lim_{t\rightarrow \infty}\|u(t)\|_{\dot{H}^{1/2}}=0.$$
Benameur \cite{JB} obtained a decay result in the Lei--Lin space. Benameur and Jlali \cite{BJ2} considered the subcritical Sobolev--Gevrey space $H^1_{a,\sigma}$, defined by the norm $\|e^{a|D|^{1/\sigma}}f\|_{H^1}$ with $\sigma>1$, and proved that
$$\limsup_{t\rightarrow \infty}\|u(t)\|_{H^1_{a,\sigma}}=0$$
for every global solution $u\in C([0,\infty),H^1_{a,\sigma})$; their proof uses the decay result of \cite{GI} together with energy estimates.
 
The purpose of the present paper is to work at the critical regularity index $s=\frac12$ with the analytic weight $e^{a|\xi|}$ (that is, $\sigma=1$). We show that every global solution $u\in C\big([0,+\infty),\dot{H}^{1/2}_{a,1}(\mathbb{R}^3)\big)$ satisfies $\|u(t)\|_{\dot{H}^{1/2}_{a,1}}\to0$ as $t\to+\infty$. This decay result is of interest for two reasons. First, since $|\xi|^{2s}\le C_{s,a}|\xi|e^{2a|\xi|}$ for $s\ge\frac12$, we have $\|u\|_{\dot{H}^{s}}\le C_{s,a}^{1/2}\|u\|_{\dot{H}^{1/2}_{a,1}}$, so that decay in $\dot{H}^{1/2}_{a,1}$ controls all derivatives and yields the rate $\|u(t)\|_{\dot{H}^s}=o\big(t^{-\frac{s-1/2}{2}}\big)$ for every $s\ge\frac12$ (Remark~\ref{appendix}), whereas decay in $\dot{H}^{1/2}$ alone, as in \cite{GI}, gives no information on derivatives. Second, working at the critical index $s=\frac12$ requires new technical ingredients. The product laws in $\dot{H}^{s}_{a,1}$ fail at the endpoint $s=\frac32$, so the classical $\dot{H}^{1/2}$ energy estimates cannot be used directly; we replace them by a product law with $s=t=1$ combined with an interpolation inequality (Lemma~\ref{lem2} and Section~4). Moreover, the passage from the decay in $\dot{H}^{1/2}$ to the decay in $\dot{H}^{1/2}_{a,1}$ relies on an analyticity estimate for small data in $\dot{H}^{1/2}$ (Lemma~\ref{theo3}), and on a uniqueness argument showing that the solution restarted at a large time belongs to the class where this estimate holds (Section~5). Along the way, we also establish a self-contained local well-posedness theory in $\dot{H}^{1/2}_{a,1}$, with a global existence criterion and a blow-up criterion (Theorems~\ref{theo1} and \ref{theo2}), and a stability result for global solutions (Theorem~\ref{theo5}).\\
We now present our main results concerning existence, blow-up criteria, asymptotic decay, and stability of solutions to the three-dimensional incompressible Navier–Stokes equations in the critical Sobolev--Gevrey space $\dot{H}^{1/2}_{a,1}(\mathbb{R}^3)$.\\
 
We begin by establishing the local well-posedness of solutions in the critical space, which forms the foundation for all subsequent results.
\begin{theorem}\label{theo1}
Let $a>0$. If $u^0\in \dot{H}^{1/2}_{a,1}$ such that ${\rm div\; }u^0=0$, then there exists a maximal time $T^*\in(0,+\infty]$ such that the Navier–Stokes system (NS) admits a unique maximal solution $$u\in C([0,T^*), \dot{H}^{1/2}_{a,1}(\mathbb{R}^3))\cap L^2([0,T^*),\dot{H}^{3/2}_{a,1}(\mathbb{R}^3)) .$$
\end{theorem}
\begin{rem}
 If $u$ is a solution of (NS) system in $C([0,T^*), \dot{H}^{1/2}_{a,1}(\mathbb{R}^3))$, then automatically $u\in  L^2_{loc}([0,T^*),\dot{H}^{3/2}_{a,1}(\mathbb{R}^3)).$
\end{rem}
 
Next, we provide a condition under which local solutions can be extended globally or may blow up in finite time.
\begin{theorem}\label{theo2}
Let $a>0$ and $u\in C([0,T^*), \dot{H}^{1/2}_{a,1}(\mathbb{R}^3))\cap L^2_{loc}([0,T^*), \dot{H}^{3/2}_{a,1}(\mathbb{R}^3))$ be the maximal solution of (NS) given by Theorem~\ref{theo1}, where $T^*\in(0,+\infty]$. Then:\\
\begin{enumerate}
 \item[(i)] If $\|u(0)\|_{\dot{H}^{1/2}_{a,1}}<\frac{1}{C_0}$, where $C_0$ is the constant of Lemma~\ref{lem3}, then $T^*=+\infty$.
 \item[(ii)] If $T^*$ is finite, then $\displaystyle\int^{T^*}_0\|u(t)\|^2_{\dot{H}^{3/2}_{a,1}}dt=+\infty.$
 \end{enumerate}
\end{theorem}
 
Having established global existence conditions, we now investigate the asymptotic behavior of global solutions.
\begin{theorem}\label{theo4}
Let $a>0$, $u\in C(\mathbb{R}^+, \dot{H}^{1/2}_{a,1}(\mathbb{R}^3))$ be a global solution of (NS), then we have:
\begin{eqnarray}\label{T4}
\lim_{t\rightarrow +\infty}\|u(t)\|_{\dot{H}^{1/2}_{a,1}}=0.
\end{eqnarray}
\end{theorem}
 \begin{rem}\label{appendix}
 The proof of this Theorem implies a result of polynomial decay in the homogeneous Sobolev spaces $\dot{H}^s(\mathbb{R}^3),$ for $s\geq \frac{1}{2},$ precisely we have:
$$
\|u(t)\|_{\dot{H}^s} =o(t^{-\frac{s-\frac{1}{2}}{2}}), \quad t\rightarrow +\infty.
$$
For the proof see Appendix.
\end{rem}
 
Finally, we address the robustness of global solutions under small perturbations of the initial data.
\begin{theorem}\label{theo5}
Let $a>0$ and let $u\in C(\mathbb{R}^+,\dot{H}^{1/2}_{a,1}(\mathbb{R}^3))$ be a global solution of (NS). Then $\Lambda:=\int^\infty_0\|u(z)\|^4_{\dot{H}^1_{a,1}}dz<\infty$, and there exist constants $K,C_1>0$, depending only on the constant $C_0$ of Lemma~\ref{lem3} ($K=216\,C_0^4$, $C_1=\frac{1}{32C_0^2}$), such that for all $v^0\in \dot{H}^{1/2}_{a,1}(\mathbb{R}^3)$ with ${\rm div}\,v^0=0$ and
$$\|v^0-u(0)\|^2_{\dot{H}^{1/2}_{a,1}}\leq C_1\,e^{-K\Lambda},$$
the Navier--Stokes system with initial data $v^0$ has a global solution $v$. Moreover, for all $t\geq 0$,
$$
\|v(t)-u(t)\|^2_{\dot{H}^{1/2}_{a,1}}+\int^t_0\|v(s)-u(s)\|^2_{\dot{H}^{3/2}_{a,1}}ds\leq \|v^0-u(0)\|^2_{\dot{H}^{1/2}_{a,1}}\,e^{K\Lambda}.
$$
\end{theorem}
This article is organized as follows. In Section 2, we give some important preliminary results. Section 3 is devoted to the existence of a local solution in the critical Sobolev--Gevrey space $\dot{H}^{1/2}_{a,1}$ (Theorem~\ref{theo1}). In Section 4, we prove the global existence criterion and the blow-up criterion (Theorem~\ref{theo2}). In Section 5, we prove the long-time decay (Theorem~\ref{theo4}). Finally, Section 6 contains the proof of the stability result for global solutions (Theorem~\ref{theo5}).
 
\section{ Notations and preliminary results}
\subsection{Notations}
In this section, we collect some notations and definitions that will be used later.\\
\begin{itemize}
\item The Fourier transformation is normalized as
$$
\mathcal{F}(f)(\xi)=\widehat{f}(\xi)=\int_{\mathbb R^3}\exp(-ix.\xi)f(x)dx,\,\,\,\xi=(\xi_1,\xi_2,\xi_3)\in\mathbb R^3.
$$
\item The inverse Fourier formula is
$$
\mathcal{F}^{-1}(g)(x)=(2\pi)^{-3}\int_{\mathbb R^3}\exp(i\xi.x)g(\xi)d\xi,\,\,\,x=(x_1,x_2,x_3)\in\mathbb R^3.
$$
\item The convolution product of a suitable pair of function $f$ and $g$ on $\mathbb R^3$ is given by
$$
(f\ast g)(x):=\int_{\mathbb R^3}f(y)g(x-y)dy.
$$
\item If $f=(f_1,f_2,f_3)$ and $g=(g_1,g_2,g_3)$ are two vector fields, we set
$$
f\otimes g:=(g_1f,g_2f,g_3f),
$$
and
$$
{\rm div}\,(f\otimes g):=({\rm div}\,(g_1f),{\rm div}\,(g_2f),{\rm div}\,(g_3f)).
$$
 
\item Let $(X,\|\cdot\|)$ be a Banach space, $1\le p\le\infty$ and $T>0$.
We denote by $L^p_T(X)$ the space of measurable functions
$f:[0,T]\to X$ such that $t\mapsto\|f(t)\|$ belongs to $L^p([0,T])$.
 
\item The homogeneous Sobolev space is
\[
\dot H^s=\{f\in\mathcal S'(\mathbb{R}^3):\ \widehat f\in L^1_{loc}(\mathbb{R}^3),\ 
|\xi|^s\widehat f\in L^2(\mathbb{R}^3)\}.
\]
 
\item For $a\ge 0$ and $s\in\mathbb{R}$, the homogeneous Sobolev--Gevrey space is
\[
\dot H^s_{a,1}(\mathbb{R}^3)=\{f\in\mathcal S'(\mathbb{R}^3):\ \widehat f\in L^1_{loc}(\mathbb{R}^3),\ 
e^{a|\xi|}|\xi|^s\widehat f\in L^2(\mathbb{R}^3)\},
\]
with the norm
\[
\|f\|_{\dot H^s_{a,1}}=\Big(\int_{\mathbb{R}^3}|\xi|^{2s}e^{2a|\xi|}|\widehat f(\xi)|^2\,d\xi\Big)^{1/2}.
\]
Equivalently, $f\in\dot H^s_{a,1}$ if and only if $e^{a|D|}f\in\dot H^s$,
where $|D|=(-\Delta)^{1/2}$. For $a=0$ we recover $\dot H^s$.
For a vector field $f=(f_1,f_2,f_3)$ we set
$\|f\|^2_{\dot H^s_{a,1}}=\sum_{i=1}^3\|f_i\|^2_{\dot H^s_{a,1}}$, and similarly for the other norms below.
 
\item We also use the following space-time spaces, whose norms are
\[
\|f\|_{\widetilde L^\infty(\dot H^{1/2})}
=\Big(\int_{\mathbb{R}^3}|\xi|\,\sup_{t\ge0}|\widehat f(t,\xi)|^2\,d\xi\Big)^{1/2},
\qquad
\|f\|_{L^2(\dot H^{3/2})}
=\Big(\int_0^\infty\!\!\int_{\mathbb{R}^3}|\xi|^{3}\,|\widehat f(t,\xi)|^2\,d\xi\,dt\Big)^{1/2}.
\]
\end{itemize}
\subsection{Preliminary results}
In this section, we recall some classical results and we give new technical lemmas.\\
It is well known that:\\
$\bullet$ The homogeneous Sobolev spaces $\dot{H}^s(\R^3)$ are Banach spaces if and only if $s<\frac{3}{2}.$\\
$\bullet$ The Sobolev--Gevrey spaces $\dot{H}^s_{a,1}(\R^3)$ are Banach spaces if and only if $s<\frac{3}{2}.$ (See \cite{BJ1})
\begin{lem}\label{lem1} (see \cite{HB})
Let $(E,\|\cdot\|)$ be a Banach space and let $B:E\times E\to E$ be a continuous bilinear map, with
$$
\|B\|=\sup_{\|x\|\le 1,\ \|y\|\le 1}\|B(x,y)\|.
$$
Let $y\in E$ be such that $4\|B\|\,\|y\|<1$. Then the equation
$$
x=y+B(x,x)
$$
has a unique solution $x$ in the closed ball $\{x\in E:\ \|x\|\le 2\|y\|\}$.
Moreover, $x$ is the only solution in the open ball $\{x\in E:\ \|x\|<\frac{1}{2\|B\|}\}$
(with the convention $\frac{1}{0}=+\infty$).
\end{lem}
\begin{lem}\label{lem2}(see \cite{HB}, \cite{JY1})
Let $a\ge 0$ and let $(s,t)\in\mathbb{R}^2$ satisfy $s<\frac32$, $t<\frac32$ and $s+t>0$.
Then there exists a constant $C=C(s,t)>0$, independent of $a$, such that for all
$u\in\dot H^{s}_{a,1}(\mathbb{R}^3)$ and $v\in\dot H^{t}_{a,1}(\mathbb{R}^3)$,
\begin{equation}\label{enq2}
\|uv\|_{\dot H^{s+t-\frac32}_{a,1}}\le C\,\|u\|_{\dot H^{s}_{a,1}}\,\|v\|_{\dot H^{t}_{a,1}} .
\end{equation}
\end{lem}
The following lemmas are inspired by \cite{FK}.
\begin{lem}\label{lem3}
Let $a\ge 0$. There exists a constant $C_0>0$ such that for all $u,v\in\dot{H}^1_{a,1}(\R^3)$ we have:
\begin{eqnarray}\label{enq3}
\|\mathbb{P}\,{\rm div}(u\otimes v)\|_{\dot{H}^{-1/2}_{a,1}}\leq C_0 \|u\|_{\dot{H}^1_{a,1}}\|v\|_{\dot{H}^1_{a,1}}.
\end{eqnarray}
\end{lem}
\begin{proof}
For $\xi\neq0$, the matrix $\delta_{jk}-\xi_j\xi_k/|\xi|^2$ is the orthogonal projection of $\mathbb{C}^3$ onto $\xi^{\perp}$,
then $$|\widehat{\mathbb{P}w}(\xi)|\le|\widehat w(\xi)| {\quad \text{and}\quad} \|\mathbb{P}w\|_{\dot H^{s}_{a,1}}\le\|w\|_{\dot H^{s}_{a,1}}.$$
Next, the $i$-th component of ${\rm div}(u\otimes v)$ is $w_i=\sum_{j}\partial_j(u_jv_i)$, it follows that 
$$|\widehat{w_i}(\xi)|\le|\xi|\sum_j|\widehat{u_jv_i}(\xi)|$$ and therefore
$$\|w_i\|_{\dot H^{-1/2}_{a,1}}\le\sum_j\|u_jv_i\|_{\dot H^{1/2}_{a,1}}.$$
Applying Lemma \ref{lem2} with $s=t=1$ (so that $s+t-\frac32=\frac12$), we obtain
\begin{align*}
\|\mathbb{P}\,{\rm div}(u\otimes v)\|_{\dot{H}^{-1/2}_{a,1}}
&\le\sum_{i,j}\|u_jv_i\|_{\dot H^{1/2}_{a,1}}\\
&\le C(1,1)\sum_{i,j}\|u_j\|_{\dot H^{1}_{a,1}}\|v_i\|_{\dot H^{1}_{a,1}}\\
&\le 3\,C(1,1)\,\|u\|_{\dot{H}^1_{a,1}}\|v\|_{\dot{H}^1_{a,1}} .
\end{align*}
Thus (\ref{enq3}) holds with $C_0=3C(1,1)$.
\end{proof}
\begin{lem}\label{lem4}
Let $a\ge0$, $T>0$, $u^0\in\dot{H}^{1/2}_{a,1}$ and $f\in L^2_T(\dot{H}^{-1/2}_{a,1})$. Let $u$ be the solution of the Cauchy problem
$$\begin{cases}
\partial_t u- \Delta u=f\\
u(0)=u^0
\end{cases}$$
given by Duhamel's formula, namely
$$
\widehat{u}(t,\xi)=e^{- t|\xi|^2}\widehat{u^0}(\xi)+\int^t_0 e^{-(t-s)|\xi|^2}\widehat{f}(s,\xi)ds .
$$
Then $u\in C([0,T], \dot{H}^{1/2}_{a,1})\cap\big(\bigcap_{2\le p\le\infty}L^p_T(\dot{H}^{\frac{1}{2}+\frac{2}{p}}_{a,1})\big)$.
Moreover, for all $t\in[0,T]$ we have the following estimates:
\begin{eqnarray}\label{enq41}
\|u(t)\|^2_{\dot{H}^{1/2}_{a,1}}+\int^t_0\|u(s)\|^2_{\dot{H}^{3/2}_{a,1}}ds\leq \|u^0\|^2_{\dot{H}^{1/2}_{a,1}}+\int^t_0\|f(s)\|^2_{\dot{H}^{-1/2}_{a,1}}ds,
\end{eqnarray}
\begin{eqnarray}\label{enq42}
\Big[\int_{\R^3} |\xi| e^{2a|\xi|}\Big(\sup_{0\leq t'\leq t}|\widehat{u}(t',\xi)|\Big)^2d\xi\Big]^{1/2}\leq \|u^0\|_{\dot{H}^{1/2}_{a,1}}+\frac{1}{\sqrt{2}}\|f\|_{L^2_t(\dot{H}^{-1/2}_{a,1})},
\end{eqnarray}
and, for every $2\le p\le\infty$,
\begin{eqnarray}\label{enq43}
\|u\|_{L_T^p(\dot{H}_{a,1}^{1/2+2/p})}\leq \|u^0\|_{\dot{H}^{1/2}_{a,1}}+\|f\|_{L^2_T(\dot{H}^{-1/2}_{a,1})}.
\end{eqnarray}
\end{lem}
\noindent Recall that $\|\nabla u\|_{\dot H^{1/2}_{a,1}}=\|u\|_{\dot H^{3/2}_{a,1}}$.
\begin{proof}
Since $f\in L^2_T(\dot{H}^{-1/2}_{a,1})$, Fubini's theorem gives
$\int_{\R^3}|\xi|^{-1}e^{2a|\xi|}\|\widehat f(\cdot,\xi)\|^2_{L^2(0,T)}d\xi<\infty$. Hence $\widehat f(\cdot,\xi)\in L^2(0,T)$ for almost every $\xi$,
and for such $\xi$ the function $t\mapsto\widehat u(t,\xi)$ is absolutely continuous.
 
\medskip\noindent{\bf Energy estimate \eqref{enq41}.} For almost every $\xi$ and almost every $t$,
$$
\frac{d}{dt}|\widehat u|^2+2|\xi|^2|\widehat u|^2=2\,{\rm Re}\big(\widehat f\,\overline{\widehat u}\big).
$$
Multiplying by $|\xi|$ and using
$2|\xi||\widehat f||\widehat u|=2\big(|\xi|^{-1/2}|\widehat f|\big)\big(|\xi|^{3/2}|\widehat u|\big)\le|\xi|^{-1}|\widehat f|^2+|\xi|^3|\widehat u|^2$,
we get
$$
\frac{d}{dt}\big(|\xi||\widehat u|^2\big)+|\xi|^3|\widehat u|^2\le|\xi|^{-1}|\widehat f|^2 .
$$
We integrate in time on $[0,t]$, multiply by $e^{2a|\xi|}$ and integrate in $\xi$. This is \eqref{enq41}.
 
\medskip\noindent{\bf Estimate \eqref{enq42}.} By Duhamel's formula and the Cauchy--Schwarz inequality in time, for every $t'\le t$,
\begin{align*}
|\widehat{u}(t',\xi)| &\leq |\widehat{u^0}(\xi)|+\int^{t'}_0 e^{-(t'-s)|\xi|^2}|\widehat{f}(s,\xi)|ds\\
                   &\leq |\widehat{u^0}(\xi)|+\Big(\int^{t'}_0 e^{-2(t'-s)|\xi|^2}ds\Big)^{\frac{1}{2}}\Big(\int^{t}_0|\widehat{f}(s,\xi)|^2ds\Big)^{\frac{1}{2}}\\
                   &\leq |\widehat{u^0}(\xi)|+\frac{1}{\sqrt{2}\,|\xi|}\|\widehat{f}(\cdot,\xi)\|_{L^2(0,t)}.
\end{align*}
Taking the supremum in $t'\in[0,t]$, multiplying by $|\xi|^{1/2}e^{a|\xi|}$ and taking the $L^2$ norm in $\xi$ (Minkowski's inequality), we obtain \eqref{enq42}.
 
\medskip\noindent{\bf Continuity in time.}
For almost every $\xi$, $t\mapsto\widehat u(t,\xi)$ is continuous, and
$|\xi|e^{2a|\xi|}|\widehat u(t,\xi)-\widehat u(t_0,\xi)|^2\le4|\xi|e^{2a|\xi|}\sup_{0\le t'\le T}|\widehat u(t',\xi)|^2$,
which is integrable in $\xi$ by \eqref{enq42}. By the dominated convergence theorem,
$u\in C([0,T];\dot H^{1/2}_{a,1})$.
 
\medskip\noindent{\bf Estimate \eqref{enq43}.}
Let $E=\|u^0\|^2_{\dot H^{1/2}_{a,1}}+\|f\|^2_{L^2_T(\dot H^{-1/2}_{a,1})}$.
By \eqref{enq41}, $\sup_{t\le T}\|u(t)\|_{\dot H^{1/2}_{a,1}}\le E^{1/2}$ and $\|u\|_{L^2_T(\dot H^{3/2}_{a,1})}\le E^{1/2}$.
This is \eqref{enq43} for $p=\infty$ and $p=2$. For $2<p<\infty$, H\"older's inequality in $\xi$ (with exponents $\frac{p}{p-2}$ and $\frac p2$)
gives the interpolation inequality
$$
\|u\|_{\dot{H}_{a,1}^{\frac{1}{2}+\frac{2}{p}}}\leq {\|u\|^{1-\frac{2}{p}}_{\dot{H}_{a,1}^{1/2}}}{\|u\|^{\frac{2}{p}}_{\dot{H}^{3/2}_{a,1}}},
$$
hence
$$
\int_0^T\|u\|^p_{\dot{H}_{a,1}^{\frac{1}{2}+\frac{2}{p}}}dt\leq {\sup_{t\le T}\|u(t)\|^{p-2}_{\dot{H}_{a,1}^{1/2}}}\int_0^T\|u\|^{2}_{\dot{H}^{3/2}_{a,1}}dt\le E^{\frac{p-2}{2}}E=E^{\frac p2}.
$$
Thus $\|u\|_{L^p_T(\dot H^{1/2+2/p}_{a,1})}\le E^{1/2}\le\|u^0\|_{\dot H^{1/2}_{a,1}}+\|f\|_{L^2_T(\dot H^{-1/2}_{a,1})}$.
\end{proof}
 
The regularizing effects of the critical space $\dot{H}^{1/2}(\mathbb{R}^3)$ of $(NS)$ equations gives us $u\in \dot{H}^{1/2}_{\sqrt{t},1}(\mathbb{R}^3):$
\begin{lem}\label{theo3}
There exist universal constants $\epsilon_0>0$, $\rho_0>0$ and $c_0>0$ such that for any initial data $u^0\in\dot{H}^{1/2}$ with ${\rm div}\,u^0=0$ and
$\|u^0\|_{\dot{H}^{1/2}}<\epsilon_0$, the Navier--Stokes system (NS) has a global solution
$u\in\widetilde{L}^\infty(\dot{H}^{1/2})\cap L^2(\dot{H}^{3/2})$, which is unique in the ball of this space with center $0$ and radius $\rho_0$, and which is analytic in the sense that:
 \begin{eqnarray}\label{T3}
  \|e^{\sqrt{t}|D|}u\|_{\widetilde{L}^\infty(\dot{H}^{1/2})}+\|e^{\sqrt{t}|D|}u\|_{L^2(\dot{H}^{3/2})}\leq c_0 \|u^0\|_{\dot{H}^{1/2}}.
  \end{eqnarray}
Here $e^{\sqrt{t}|D|}u$ denotes the function whose Fourier transform is $e^{\sqrt{t}|\xi|}\widehat{u}(t,\xi)$.
\end{lem}
\noindent{\bf Proof.}
This proof is inspired by the work of Bae \cite{BAE} on the Lei--Lin space; we give a self-contained argument in $\dot{H}^{1/2}$. It is done in three steps.\\
We first take the Fourier transform to the integral form of Navier-Stokes equation:
\begin{eqnarray}\label{enqT40}
\widehat{u}(t,\xi)=e^{- t|\xi|^2}\widehat{u}^0-\int^t_0 e^{- (t-s)|\xi|^2}\widehat{f}(s,\xi)ds.
\end{eqnarray}
{\bf Notation.} We set $X:=\widetilde{L}^\infty(\dot{H}^{1/2})\cap L^2(\dot{H}^{3/2})$ (time interval $[0,\infty)$), with the norm
$$
\|u\|_X:=\|u\|_{\widetilde{L}^\infty(\dot{H}^{1/2})}+\|u\|_{L^2(\dot{H}^{3/2})}.
$$
$X$ is a Banach space: it is the intersection of two weighted mixed-norm spaces of functions $\widehat u(t,\xi)$, and $\widehat u(t,\cdot)\in L^1_{loc}$ by the first factor. Clearly $\sup_{t\ge0}\|u(t)\|_{\dot{H}^{1/2}}\le\|u\|_{\widetilde{L}^\infty(\dot{H}^{1/2})}$.
In (\ref{enqT40}), $\widehat f=\widehat{\mathbb{P}\,{\rm div}(u\otimes u)}$, so that $|\widehat f(s,\xi)|\le|\xi|\,|\widehat{u\otimes u}(s,\xi)|$. We denote by $C_{s,t}$ the constant of Lemma \ref{lem2} (with $a=0$).\\
 
{\bf Step 1:} we start by estimating $u$ in $\widetilde{L}^\infty(\dot{H}^{1/2}).$ Multiplying (\ref{enqT40}) by $|\xi|^{1/2}$, we get:
\begin{align*}
|\xi|^{1/2}|\widehat{u}(t,\xi)|&\leq |\xi|^{1/2}|\widehat{u}^0(\xi)|+\int^t_0e^{-(t-s)|\xi|^2}|\xi|^{3/2}|\widehat{u\otimes u}(s,\xi)|ds\\
                    &\leq |\xi|^{1/2}|\widehat{u}^0(\xi)|+\int^t_0|\xi|^2e^{-(t-s)|\xi|^2}ds\ \sup_{0\leq s<\infty}|\xi|^{-1/2}|\widehat{u\otimes u}(s,\xi)|\\
                    &\leq |\xi|^{1/2}|\widehat{u}^0|+ |\xi|^{-1/2}\sup_{0\leq s<\infty}|\widehat{u\otimes u}(s,\xi)|.
\end{align*}
The right-hand side does not depend on $t$, hence the same bound holds for $|\xi|^{1/2}\sup_{t\ge0}|\widehat u(t,\xi)|$. Let $U(\xi)=\sup_{t\ge0}|\widehat u(t,\xi)|$. Then for all $s$ $$|\widehat{u\otimes u}(s,\xi)|\le(2\pi)^{-3}(U*U)(\xi),$$  and $U$ is the Fourier transform of a tempered distribution $\widetilde U$ with $$\|\widetilde U\|_{\dot H^{1/2}}=\|u\|_{\widetilde{L}^\infty(\dot{H}^{1/2})}$$ (as in the proof of Lemma \ref{lem2}). Taking the $L^2$ norm with respect to the frequency variable $\xi$ and using Lemma \ref{lem2} with $s=t=\frac12$, we obtain:
\begin{eqnarray}\label{enqT41}
\|u\|_{\widetilde{L}^\infty(\dot{H}^{1/2})}\leq \|u^0\|_{\dot{H}^{1/2}}+C_{\frac{1}{2},\frac{1}{2}} \|u\|^2_{\widetilde{L}^\infty(\dot{H}^{1/2})}.
\end{eqnarray}
Now, we estimate $u$ in $L^2( \dot{H}^{3/2})$. Multiplying (\ref{enqT40}) by $|\xi|^{3/2}$, we get:
$$
|\xi|^{3/2}|\widehat{u}(t,\xi)|\leq |\xi|e^{-t|\xi|^2}|\xi|^{1/2}|\widehat{u}^0(\xi)|+\int^t_0|\xi|^2e^{-(t-s)|\xi|^2}\,|\xi|^{1/2}|\widehat{u\otimes u}(s,\xi)|ds .
$$
Taking the $L^2$ norm with respect to time and using Young's inequality, we deduce:
\begin{align*}
\Big(\int^\infty_0|\xi|^{3}|\widehat{u}(t,\xi)|^2dt\Big)^{1/2}&\leq \Big(\int^\infty_0 |\xi|^2e^{-2s|\xi|^2}ds\Big)^{1/2}|\xi|^{1/2}|\widehat{u}^0(\xi)|\\
&\quad+ \int^\infty_0|\xi|^2e^{-s|\xi|^2}ds\Big(\int^\infty_0|\xi||\widehat{u\otimes u}(s,\xi)|^2ds\Big)^{1/2}\\
                                              &\leq  |\xi|^{1/2}|\widehat{u}^0(\xi)|+\Big(\int^\infty_0|\xi||\widehat{u\otimes u}(s,\xi)|^2ds\Big)^{1/2}.
\end{align*}
Taking the $L^2$ norm in $\xi$, using Lemma \ref{lem2} with $s=t=1$ (so that $s+t-\frac32=\frac12$) and the inequality $\|w\|^2_{\dot{H}^{1}}\le\|w\|_{\dot{H}^{1/2}}\|w\|_{\dot{H}^{3/2}}$ (Cauchy--Schwarz in $\xi$), we obtain:
\begin{align*}
\|u\|_{L^2(\dot{H}^{3/2})}&\leq  \|u^0\|_{\dot{H}^{1/2}}+\Big(\int^\infty_0\int_{\mathbb{R}^3} |\xi||\widehat{u\otimes u}(s,\xi)|^2d\xi ds\Big)^{1/2}\\
                                   &=  \|u^0\|_{\dot{H}^{1/2}}+\Big(\int^\infty_0\|(u\otimes u)(s)\|^2_{\dot{H}^{1/2}}ds\Big)^{1/2}\\
                                   &\leq \|u^0\|_{\dot{H}^{1/2}}+C_{1,1}\Big(\int^\infty_0\|u(s)\|^4_{\dot{H}^{1}}ds\Big)^{1/2}\\
                                   &\leq \|u^0\|_{\dot{H}^{1/2}}+C_{1,1}\Big(\int^\infty_0\|u(s)\|^2_{\dot{H}^{1/2}}\|u(s)\|^2_{\dot{H}^{3/2}}ds\Big)^{1/2}\\
                                   &\leq \|u^0\|_{\dot{H}^{1/2}}+C_{1,1}\Big(\sup_{t\ge0}\int_{\mathbb{R}^3} |\xi||\widehat{u}(t,\xi)|^2d\xi\ \int^\infty_0\|u(s)\|^2_{\dot{H}^{3/2}}ds\Big)^{1/2}\\
                                   &\leq \|u^0\|_{\dot{H}^{1/2}}+C_{1,1}\|u\|_{L^2(\dot{H}^{3/2})}\Big(\int_{\mathbb{R}^3} |\xi|\sup_{t\ge0}|\widehat{u}(t,\xi)|^2d\xi\Big)^{1/2},
\end{align*}
which yields,
 \begin{eqnarray}\label{enqT42}
 \|u\|_{L^2(\dot{H}^{3/2})}\leq \|u^0\|_{\dot{H}^{1/2}}+C_{1,1}\|u\|_{L^2( \dot{H}^{3/2})}\|u\|_{\widetilde{L}^\infty(\dot{H}^{1/2})}.
\end{eqnarray}
The same computations, with $\widehat{u\otimes u}$ replaced by $\widehat{u\otimes v}$ (and the Cauchy--Schwarz inequality in time in the last lines), give for the bilinear operator $B(u,v)=-\int_0^te^{(t-s)\Delta}\mathbb{P}\,{\rm div}(u\otimes v)(s)\,ds$ and all $u,v\in X$:
\begin{eqnarray}\label{enqT45}
\|B(u,v)\|_{\widetilde{L}^\infty(\dot{H}^{1/2})}\le C_{\frac12,\frac12}\|u\|_{\widetilde{L}^\infty(\dot{H}^{1/2})}\|v\|_{\widetilde{L}^\infty(\dot{H}^{1/2})},\\
\nonumber\|B(u,v)\|_{L^2(\dot{H}^{3/2})}\le C_{1,1}\|u\|_X\|v\|_X,
\end{eqnarray}
hence $\|B(u,v)\|_X\le C\|u\|_X\|v\|_X$, where $C=C_{\frac12,\frac12}+C_{1,1}$.
 
{\bf Step 2:} Combining (\ref{enqT41}) and (\ref{enqT42}), we get:
\begin{eqnarray}\label{enqT43}
\|u\|_{\widetilde{L}^\infty(\dot{H}^{1/2})}+\|u\|_{L^2(\dot{H}^{3/2})}\leq 2\|u^0\|_{\dot{H}^{1/2}}+C(\|u\|_{\widetilde{L}^\infty(\dot{H}^{1/2})}+\|u\|_{L^2(\dot{H}^{3/2})})^2,
\end{eqnarray}
with $C=C_{\frac12,\frac12}+C_{1,1}.$\\
Let $r=\frac{1}{4C}$, and let $\|u^0\|_{\dot{H}^{1/2}}<\epsilon_1:=\frac{3r}{8}$. We take:
$$
B_r=\{u\in \widetilde{L}^\infty([0,\infty[,\dot{H}^{1/2})\cap L^2([0,\infty[, \dot{H}^{3/2})\ /\ \|u\|_{\widetilde{L}^\infty(\dot{H}^{1/2})}+\|u\|_{L^2(\dot{H}^{3/2})}\leq r\}.
$$
We consider the application $\psi$ defined by:
$$
\psi(u)=e^{t\Delta}u^0-\int^t_0 e^{(t-s)\Delta} \mathbb{P}({\rm div}(u\otimes u)(s))ds.
$$
Since $\|e^{t\Delta}u^0\|_{\widetilde{L}^\infty(\dot{H}^{1/2})}=\|u^0\|_{\dot{H}^{1/2}}$ and $\|e^{t\Delta}u^0\|_{L^2(\dot{H}^{3/2})}=\frac{1}{\sqrt2}\|u^0\|_{\dot{H}^{1/2}},$ we have $$\|e^{t\Delta}u^0\|_X\le2\|u^0\|_{\dot{H}^{1/2}}.$$ Then, by (\ref{enqT45}), we have:
$$
\|\psi(u)\|_{\widetilde{L}^\infty(\dot{H}^{1/2})}+\|\psi(u)\|_{L^2(\dot{H}^{3/2})}\leq  2\|u^0\|_{\dot{H}^{1/2}}+C (\|u\|_{\widetilde{L}^\infty(\dot{H}^{1/2})}+\|u\|_{L^2(\dot{H}^{3/2})})^2,
$$
which yields, for $u\in B_r$,
$$
\|\psi(u)\|_{\widetilde{L}^\infty(\dot{H}^{1/2})}+\|\psi(u)\|_{L^2(\dot{H}^{3/2})} <2\epsilon_1+Cr^2=\frac{3r}{4}+\frac{r}{4}=r,
$$
then $$\psi(B_r)\subset B_r.$$
We have, for all $u_1,u_2 \in B_r$:
\begin{align*}
\|\psi(u_1)-\psi(u_2)\|_{X}&\leq \|B(u_1-u_2,u_1)+B(u_2,u_1-u_2)\|_{X}\\
                                                     &\leq C (\|u_1\|_{X}+\|u_2\|_{X}) \|u_1-u_2\|_{X}\\
                                                     &\leq 2Cr \|u_1-u_2\|_{X}\\
                                                     &\leq \frac{1}{2} \|u_1-u_2\|_{X}.
\end{align*}
Which implies the existence of a unique global solution $u\in B_r\subset \widetilde{L}^\infty(\dot{H}^{1/2})\cap L^2(\dot{H}^{3/2})$ for small initial data in $\dot{H}^{1/2}(\mathbb{R}^3)$, and, by (\ref{enqT43}) and $\|u\|_X\le r$, we get $\|u\|_X\le2\|u^0\|_{\dot{H}^{1/2}}+\frac14\|u\|_X$, that is
$$\|u\|_{\widetilde{L}^\infty(\dot{H}^{1/2})}+\|u\|_{L^2( \dot{H}^{3/2})}\leq \frac83\|u^0\|_{\dot{H}^{1/2}}.$$
 
{\bf Step 3:} Let $Y$ be the space of functions $u$ such that $V:=e^{\sqrt{t}|D|}u\in X$, with the norm $\|u\|_Y:=\|V\|_X$, where $\widehat V(t,\xi)=e^{\sqrt{t}|\xi|}\widehat u(t,\xi)$. The map $u\mapsto V$ is an isometry of $Y$ onto $X$, so $Y$ is a Banach space, and $\|u\|_X\le\|u\|_Y$.
Multiplying (\ref{enqT40}) by $e^{\sqrt{t}|\xi|}$ we obtain:
\begin{align*}
e^{\sqrt{t}|\xi|}|\widehat{u}(t,\xi)| &\leq e^{\sqrt{t}|\xi|-t|\xi|^2}|\widehat{u}^0|+ \int^t_0 e^{(\sqrt{t}-\sqrt{s})|\xi|-(t-s)|\xi|^2}\,|\xi|\,e^{\sqrt{s}|\xi|}|\widehat{u\otimes u}(s,\xi)|ds.
\end{align*}
Since $\sqrt t-\sqrt s\le\sqrt{t-s}$ and $x-x^2\le\frac12-\frac{x^2}{2}$ for all $x\in\mathbb{R}$, we have $$e^{(\sqrt t-\sqrt s)|\xi|-(t-s)|\xi|^2}\le\sqrt{e}\,e^{-\frac12(t-s)|\xi|^2} {\quad \text{and}\quad} e^{\sqrt t|\xi|-t|\xi|^2}\le\sqrt e\,e^{-\frac12t|\xi|^2}.$$ Moreover, since $|\xi|\le|\xi-\eta|+|\eta|$,
\begin{align*}
e^{\sqrt{s}|\xi|}|\widehat{u\otimes u}(s,\xi)|&\le(2\pi)^{-3}\int e^{\sqrt{s}|\xi-\eta|}|\widehat{u}(s,\xi-\eta)|\,e^{\sqrt{s}|\eta|}|\widehat{u}(s,\eta)|\,d\eta\\&=(2\pi)^{-3}\big(|\widehat{V}(s,\cdot)|*|\widehat{V}(s,\cdot)|\big)(\xi),
\end{align*}
with $V(t,.)= e^{\sqrt{t}|D|}u(t,.)$. Then we get:
\begin{align*}
|\widehat{V}(t,\xi)|&\leq \sqrt{e}\Big(e^{-\frac{1}{2}t|\xi|^2}|\widehat{u}^0|+\int^t_0|\xi|e^{-\frac{1}{2}(t-s)|\xi|^2}\,(2\pi)^{-3}\big(|\widehat{V}(s,\cdot)|*|\widehat{V}(s,\cdot)|\big)(\xi)\,ds\Big),
\end{align*}
which is the situation of Step 1 with the kernel $e^{-\tau|\xi|^2}$ replaced by $e^{-\frac12\tau|\xi|^2}$ and the factor $\sqrt e$. Since $\int_0^\infty e^{-\frac12\tau|\xi|^2}d\tau=2|\xi|^{-2}$ and $\int_0^\infty|\xi|^2e^{-\frac12\tau|\xi|^2}d\tau=2$, the constants $1$ of Step 1 become $2$, and the same computations (for the bilinear operator, with $V_u=e^{\sqrt t|D|}u$ and $V_v=e^{\sqrt t|D|}v$ in the place of $u$ and $v$) give, for all $u,v\in Y$,
\begin{eqnarray}\label{enqT46}
\|B(u,v)\|_{Y}\le C_Y\|u\|_Y\|v\|_Y,\qquad C_Y:=2\sqrt{e}\,C.
\end{eqnarray}
For the free evolution, $\widehat{e^{\sqrt t|D|}e^{t\Delta}u^0}(t,\xi)=e^{\sqrt t|\xi|-t|\xi|^2}\widehat{u^0}(\xi)$. Since $x-x^2\le\frac14$, we get
$$
\|e^{\sqrt t|D|}e^{t\Delta}u^0\|_{\widetilde{L}^\infty(\dot{H}^{1/2})}\le e^{1/4}\|u^0\|_{\dot{H}^{1/2}}.
$$
With the change of variable $\tau=t|\xi|^2$ and the inequality $2\sqrt\tau-2\tau\le1-\tau$,
\begin{align*}
\int_0^\infty e^{2\sqrt t|\xi|-2t|\xi|^2}dt&=|\xi|^{-2}\int_0^\infty e^{2\sqrt\tau-2\tau}d\tau\\&\le e\,|\xi|^{-2},
\end{align*}
hence $\|e^{\sqrt t|D|}e^{t\Delta}u^0\|_{L^2(\dot{H}^{3/2})}\le e^{1/2}\|u^0\|_{\dot{H}^{1/2}}$. Consequently
\begin{eqnarray}\label{enqT47}
\|e^{t\Delta}u^0\|_Y<3\|u^0\|_{\dot{H}^{1/2}}.
\end{eqnarray}
Let $r_Y=\frac{1}{4C_Y}$, $\epsilon_0:=\frac{r_Y}{4}=\frac{1}{32\sqrt e\,C}$, and $\|u^0\|_{\dot{H}^{1/2}}<\epsilon_0$. Let $B^Y_{r_Y}=\{u\in Y:\ \|u\|_Y\le r_Y\}$ and let $\psi$ be as in Step 2. By (\ref{enqT46}) and (\ref{enqT47}), for $u\in B^Y_{r_Y}$,
\begin{align*}
\|\psi(u)\|_Y&\le3\|u^0\|_{\dot{H}^{1/2}}+C_Y\|u\|_Y^2\\&<\frac{3r_Y}{4}+\frac{r_Y}{4}\\&=r_Y,
\end{align*}
and, for $u_1,u_2\in B^Y_{r_Y}$, \begin{align*}
\|\psi(u_1)-\psi(u_2)\|_Y&\le2C_Yr_Y\|u_1-u_2\|_Y\\&=\frac12\|u_1-u_2\|_Y.\end{align*}
Hence $\psi$ has a unique fixed point $u\in B^Y_{r_Y}$, and from \begin{align*} \|u\|_Y&\le3\|u^0\|_{\dot{H}^{1/2}}+C_Yr_Y\|u\|_Y\\&=3\|u^0\|_{\dot{H}^{1/2}}+\frac14\|u\|_Y
\end{align*}
we get
\begin{eqnarray}\label{enqT44}
\|e^{\sqrt{t}|D|}u\|_{\widetilde{L}^\infty(\dot{H}^{1/2})}+\|e^{\sqrt{t}|D|}u\|_{L^2(\dot{H}^{3/2})}=\|u\|_Y\leq 4\|u^0\|_{\dot{H}^{1/2}}.
\end{eqnarray}
It remains to check that $u$ is the solution of Step 2. Since $\epsilon_0<\epsilon_1$, Step 2 applies; moreover $\|u\|_X\le\|u\|_Y\le r_Y=\frac{1}{8\sqrt e\,C}<r,$ then $u\in B_r$, where the solution of Step 2 is unique. This proves the lemma, with $\rho_0=r=\frac{1}{4C}$ and $c_0=4$.
 
\section{ Proof of Theorem~\ref{theo1}}
We follow the strategy of Fujita and Kato \cite{FK}, who proved the existence of Navier--Stokes solutions in the critical space $\dot{H}^{1/2}$.
By \eqref{NS-integral}, a function $u$ solves (NS) in the mild sense if and only if $u=e^{t\Delta}u^0+B(u,u)$, where $B(u,v)$ is the solution of the heat equation
\begin{eqnarray}\label{enqT1}
\begin{cases}
\partial_t B-\Delta B=-\mathbb{P}\,{\rm div}(u\otimes v),\\
B(0)=0.
\end{cases}
\end{eqnarray}
Let $a>0$, $T>0$ and $E_T:=L^4_T(\dot{H}^1_{a,1})$.
 
\medskip\noindent{\bf Step 1: the bilinear operator.}
By estimate \eqref{enq43} of Lemma \ref{lem4} with $p=4$, $u^0=0$ and $f=-\mathbb{P}\,{\rm div}(u\otimes v)$, then Lemma \ref{lem3} and the Cauchy--Schwarz inequality in time, we get
\begin{align*}
\|B(u,v)\|_{E_T}&\leq \|\mathbb{P}\,{\rm div}(u\otimes v)\|_{L^2_T(\dot{H}^{-1/2}_{a,1})}\\
&\leq C_0\Big(\int^T_0 \|u(s)\|^2_{\dot{H}^1_{a,1}}\|v(s)\|^2_{\dot{H}^1_{a,1}}ds\Big)^{1/2}\\
&\leq C_0 \|u\|_{E_T}\|v\|_{E_T}.
\end{align*}
Hence $B$ is a continuous bilinear map on $E_T$ with $\|B\|\le C_0$, uniformly in $T$.
 
\medskip\noindent{\bf Step 2: the free evolution.}
For every $T>0$,
\begin{eqnarray}\label{enqT2}
\|e^{ t\Delta}u^0\|_{E_T}\leq \frac{1}{\sqrt2}\|u^0\|_{\dot{H}^{1/2}_{a,1}}\leq \|u^0\|_{\dot{H}^{1/2}_{a,1}}.
\end{eqnarray}
Indeed, by Minkowski's inequality in $L^2_t$,
\begin{align*}
\|e^{ t\Delta}u^0\|^2_{E_T}
&=\Big\|\int_{\R^3}|\xi|^2e^{-2t|\xi|^2}e^{2a|\xi|}|\widehat{u^0}(\xi)|^2d\xi\Big\|_{L^2(0,T)}\\
&\le\int_{\R^3}\big\||\xi|^2e^{-2t|\xi|^2}\big\|_{L^2(0,\infty)}e^{2a|\xi|}|\widehat{u^0}(\xi)|^2d\xi\\
&=\frac12\|u^0\|^2_{\dot{H}^{1/2}_{a,1}},
\end{align*}
because $\int_0^\infty|\xi|^4e^{-4t|\xi|^2}dt=\frac{|\xi|^2}{4}$.
 
\medskip\noindent{\bf Step 3: small data.}
If $\|u^0\|_{\dot{H}^{1/2}_{a,1}}<\frac{1}{4C_0}$, then by \eqref{enqT2}, $\|e^{ t\Delta}u^0\|_{E_T}<\frac{1}{4C_0}\le\frac{1}{4\|B\|}$ for every $T>0$.
According to Lemma \ref{lem1}, there exists a unique solution $u$ of $u=e^{t\Delta}u^0+B(u,u)$ in the ball of $E_T$ with center $0$ and radius $\frac{1}{2C_0}$.
 
\medskip\noindent{\bf Step 4: large data.}
Let $u^0\in\dot{H}^{1/2}_{a,1}$ with $u^0\neq0$ and let $\varepsilon=\frac{1}{8C_0}$. Since $|\xi|^{1/2}e^{a|\xi|}\widehat{u^0}\in L^2(\R^3)$, there exists $\rho=\rho_{u^0}>0$ such that
$$
\Big(\int_{|\xi|>\rho}|\xi|e^{2a|\xi|}|\widehat{u^0}(\xi)|^2d\xi\Big)^{1/2}< \varepsilon .
$$
Write $u^0=u^0_{h}+u^0_{\ell}$ with $\widehat{u^0_{\ell}}=\mathbf 1_{\{|\xi|\le\rho\}}\widehat{u^0}$.
By \eqref{enqT2} applied to $u^0_h$, we have $\|e^{ t\Delta}u^0_h\|_{E_T}<\varepsilon$ for every $T$.
For the low frequencies, $|\xi|^2e^{-2t|\xi|^2}\le\rho\,|\xi|$ when $|\xi|\le\rho$, so that
$\|e^{ t\Delta}u^0_{\ell}\|^2_{\dot{H}^{1}_{a,1}}\le\rho\,\|u^0\|^2_{\dot{H}^{1/2}_{a,1}}$ for every $t\ge0$, and therefore
$$
\|e^{t\Delta}u^0_{\ell}\|_{E_T}\leq T^{\frac14}\rho^{\frac12}\|u^0\|_{\dot{H}^{1/2}_{a,1}}.
$$
Thus, if
\begin{eqnarray}\label{enqT3}
 T\leq \Big(\frac{\varepsilon}{\rho^{1/2}_{u^0}\|u^0\|_{\dot{H}^{1/2}_{a,1}}}\Big)^4,
\end{eqnarray}
then $\|e^{t\Delta}u^0\|_{E_T}<2\varepsilon=\frac{1}{4C_0}$, and Lemma \ref{lem1} gives, as in Step 3, a unique solution of $u=e^{t\Delta}u^0+B(u,u)$ in the ball of $E_T$ with center $0$ and radius $\frac{1}{2C_0}$.
 
\medskip\noindent{\bf Step 5: regularity.}
Let $u\in E_T$ be a solution given by Step 3 or Step 4. By Lemma \ref{lem3}, $\mathbb{P}\,{\rm div}(u\otimes u)\in L^2_T(\dot{H}^{-1/2}_{a,1})$. Since $u$ solves the heat equation with data $u^0\in\dot{H}^{1/2}_{a,1}$ and right-hand side $-\mathbb{P}\,{\rm div}(u\otimes u)$, Lemma \ref{lem4} implies that
$$
 u\in C([0,T],\dot{H}^{1/2}_{a,1})\cap L^2_T(\dot{H}^{3/2}_{a,1}).
$$
Moreover, ${\rm div}\,u=0$ because ${\rm div}\,u^0=0$ and $\mathbb{P}$ takes values in divergence-free fields.
 
\section{ Proof of Theorem~\ref{theo2}}
Let $u\in C([0,T^*),\dot{H}^{1/2}_{a,1})\cap L^2_{loc}([0,T^*),\dot{H}^{3/2}_{a,1})$ be the maximal solution of (NS). Recall that $E_T:=L^4_T(\dot{H}^1_{a,1})$ (Section 3).
 
\medskip\noindent{\bf Two elementary inequalities.} For $w\in\dot{H}^{1/2}_{a,1}\cap\dot{H}^{3/2}_{a,1}$, the Cauchy--Schwarz inequality in $\xi$ gives
\begin{eqnarray}\label{enqS0}
\|w\|^2_{\dot{H}^1_{a,1}}\le\|w\|_{\dot{H}^{1/2}_{a,1}}\|w\|_{\dot{H}^{3/2}_{a,1}} .
\end{eqnarray}
Moreover, by Lemma \ref{lem2} with $s=t=1$ (so that $s+t-\frac32=\frac12$), exactly as in the proof of Lemma \ref{lem3}, we have $$\|u\otimes u\|_{\dot{H}^{1/2}_{a,1}}\le C_0\|u\|^2_{\dot{H}^1_{a,1}},$$ with $C_0=3C(1,1)$. Hence, by (\ref{enqS0}),
\begin{eqnarray}\label{enqS0b}
\|u\otimes u\|_{\dot{H}^{1/2}_{a,1}}\le C_0\|u\|_{\dot{H}^{1/2}_{a,1}}\|u\|_{\dot{H}^{3/2}_{a,1}}.
\end{eqnarray}
 
\medskip\noindent{\bf Energy inequality.} For $t<T^*$, (\ref{enqS0}) gives
$$
\int_0^t\|u(s)\|^4_{\dot{H}^1_{a,1}}ds\le\sup_{s\le t}\|u(s)\|^2_{\dot{H}^{1/2}_{a,1}}\int_0^t\|u(s)\|^2_{\dot{H}^{3/2}_{a,1}}ds<\infty,
$$
so that $\mathbb{P}\,{\rm div}(u\otimes u)\in L^2_{loc}([0,T^*),\dot{H}^{-1/2}_{a,1})$ by (\ref{enq3}). As in the proof of (\ref{enq41}) in Lemma \ref{lem4}, we use that $t\mapsto\widehat u(t,\xi)$ is absolutely continuous and satisfies $$\frac{d}{dt}|\widehat u|^2+2|\xi|^2|\widehat u|^2=-2\,{\rm Re}\big(\widehat{\mathbb{P}\,{\rm div}(u\otimes u)}\ \overline{\widehat u}\big).$$ We multiply by $|\xi|e^{2a|\xi|}$ and integrate in $\xi$ and $t$. Since $\mathbb{P}$ is self-adjoint and $\mathbb{P}u=u$, we obtain, for $0\le t_1\le t_2<T^*$,
\begin{eqnarray}\label{enqS1}
\|u(t_2)\|^2_{\dot{H}^{1/2}_{a,1}}+2\int_{t_1}^{t_2}\|u(s)\|^2_{\dot{H}^{3/2}_{a,1}}ds=\|u(t_1)\|^2_{\dot{H}^{1/2}_{a,1}}-2\int_{t_1}^{t_2}\langle{\rm div}(u\otimes u),u\rangle_{\dot{H}^{1/2}_{a,1}}ds,
\end{eqnarray}
where $\langle f,g\rangle_{\dot{H}^{1/2}_{a,1}}={\rm Re}\int_{\mathbb{R}^3}|\xi|e^{2a|\xi|}\widehat f\cdot\overline{\widehat g}\,d\xi$. Moreover, by the Cauchy--Schwarz inequality and (\ref{enqS0b}),
\begin{align*}
|\langle{\rm div}(u\otimes u),u\rangle_{\dot{H}^{1/2}_{a,1}}|&\le\int_{\mathbb{R}^3}|\xi|^{1/2}e^{a|\xi|}|\widehat{u\otimes u}|\ |\xi|^{3/2}e^{a|\xi|}|\widehat u|\,d\xi\\
&\le\|u\otimes u\|_{\dot{H}^{1/2}_{a,1}}\|u\|_{\dot{H}^{3/2}_{a,1}}\le C_0\|u\|_{\dot{H}^{1/2}_{a,1}}\|u\|^2_{\dot{H}^{3/2}_{a,1}} .
\end{align*}
Therefore, for $0\le t_1\le t_2<T^*$,
\begin{eqnarray}\label{enqS2}
\|u(t_2)\|^2_{\dot{H}^{1/2}_{a,1}}+2\int_{t_1}^{t_2}\|u(s)\|^2_{\dot{H}^{3/2}_{a,1}}ds\leq \|u(t_1)\|^2_{\dot{H}^{1/2}_{a,1}}+2C_0\int_{t_1}^{t_2}\|u(s)\|_{\dot{H}^{1/2}_{a,1}}\|u(s)\|^2_{\dot{H}^{3/2}_{a,1}}ds.
\end{eqnarray}
 
\medskip\noindent{\bf Proof of (ii).} We argue by contradiction. Assume that $T^*<\infty$ and $\int_0^{T^*}\|u(t)\|^2_{\dot{H}^{3/2}_{a,1}}dt<\infty$.
 
\smallskip\noindent{\it Step 1: $u$ is bounded in $\dot{H}^{1/2}_{a,1}$.} Let $T\in(0,T^*)$ be such that $\int_T^{T^*}\|u(t)\|^2_{\dot{H}^{3/2}_{a,1}}dt\le\frac{1}{4C_0}$. For $T\le z\le t<T^*$, (\ref{enqS2}) gives
\begin{align*}
\|u(z)\|^2_{\dot{H}^{1/2}_{a,1}}+2\int_T^z\|u(s)\|^2_{\dot{H}^{3/2}_{a,1}}ds&\le\|u(T)\|^2_{\dot{H}^{1/2}_{a,1}}+2C_0\sup_{T\le s\le t}\|u(s)\|_{\dot{H}^{1/2}_{a,1}}\int_T^{T^*}\|u(s)\|^2_{\dot{H}^{3/2}_{a,1}}ds\\
&\le\|u(T)\|^2_{\dot{H}^{1/2}_{a,1}}+\frac12\sup_{T\le s\le t}\|u(s)\|_{\dot{H}^{1/2}_{a,1}}.
\end{align*}
Hence $S_t:=\sup_{T\le s\le t}\|u(s)\|_{\dot{H}^{1/2}_{a,1}}$ satisfies $S_t^2\le\|u(T)\|^2_{\dot{H}^{1/2}_{a,1}}+\frac12S_t$, that is
$$
S_t\le C_T:=\frac14+\sqrt{\frac{1}{16}+\|u(T)\|^2_{\dot{H}^{1/2}_{a,1}}}.
$$
Thus $\|u(t)\|_{\dot{H}^{1/2}_{a,1}}\le M:=\max\big(\sup_{[0,T]}\|u\|_{\dot{H}^{1/2}_{a,1}},\,C_T\big)$ for all $t\in[0,T^*)$.
 
\smallskip\noindent{\it Step 2: $u\in L^4([0,T^*),\dot{H}^1_{a,1})$.} By (\ref{enqS0}),
$$
\int_0^{T^*}\|u(t)\|^4_{\dot{H}^1_{a,1}}dt\le M^2\int_0^{T^*}\|u(t)\|^2_{\dot{H}^{3/2}_{a,1}}dt<\infty .
$$
Hence there exists $t_0\in(0,T^*)$ such that $\|u\|_{L^4([t_0,T^*),\dot{H}^1_{a,1})}\le\varepsilon:=\frac{1}{8C_0}$.
 
\smallskip\noindent{\it Step 3: extension beyond $T^*$.} Let $T_1=T^*-t_0$, $w_0=u(t_0)\in\dot{H}^{1/2}_{a,1}$, $y(t)=e^{t\Delta}w_0$ and $\widetilde u(t)=u(t_0+t)$ for $0\le t<T_1$. By the semigroup property, $\widetilde u=y+B(\widetilde u,\widetilde u)$ on $[0,T_1)$. By Step 1 of Section 3, $\|B(\widetilde u,\widetilde u)\|_{E_{T_1}}\le C_0\|\widetilde u\|^2_{E_{T_1}}$, hence
\begin{align*}
\|y\|_{E_{T_1}}&\le\|\widetilde u\|_{E_{T_1}}+C_0\|\widetilde u\|^2_{E_{T_1}}\\&\le\varepsilon+C_0\varepsilon^2\\&=\frac{9}{64C_0}.
\end{align*}
On the other hand, the bound (\ref{enqT2}) is uniform in $T$, so $y\in L^4(0,\infty;\dot{H}^1_{a,1})$. Hence there is $\delta>0$ such that $\|y\|_{L^4(T_1,T_1+\delta;\dot{H}^1_{a,1})}\le\frac{3}{64C_0}$. With $T_2:=T_1+\delta$ we get
$$
\|y\|_{E_{T_2}}\le\frac{1}{4\|B\|}.
$$
By Lemma \ref{lem1} in $E_{T_2}$, there is a unique $w\in E_{T_2}$ such that $w=y+B(w,w)$ and $\|w\|_{E_{T_2}}\le2\|y\|_{E_{T_2}}<\frac{1}{2C_0}$. The restriction of $w$ to $[0,T_1)$ and $\widetilde u$ both solve this equation in $E_{T_1}$ and have norm smaller than $\frac{1}{2C_0}$; by the uniqueness part of Lemma \ref{lem1} (in $E_{T_1}$, where $\|B\|\le C_0$), they coincide. By Step 5 of Section 3, $w\in C([0,T_2],\dot{H}^{1/2}_{a,1})\cap L^2_{T_2}(\dot{H}^{3/2}_{a,1})$, and $w(0)=u(t_0)$. Therefore the function equal to $u$ on $[0,t_0]$ and to $w(\cdot-t_0)$ on $[t_0,T^*+\delta)$ is a solution of (NS) on $[0,T^*+\delta)$ in the same class, which coincides with $u$ on $[0,T^*)$. This contradicts the maximality of $T^*$, and proves (ii).
 
\medskip\noindent{\bf Proof of (i).} Assume that $\|u(0)\|_{\dot{H}^{1/2}_{a,1}}<\frac{1}{C_0}$ and set $\kappa:=1-C_0\|u(0)\|_{\dot{H}^{1/2}_{a,1}}>0$. Let
$$
T_3=\sup\big\{t\in[0,T^*):\ \|u(s)\|_{\dot{H}^{1/2}_{a,1}}<\tfrac{1}{C_0}\ \text{for all } s\le t\big\},
$$
which is positive by continuity. By (\ref{enqS2}), for $0\le t_1\le t_2<T_3$,
$$
\|u(t_2)\|^2_{\dot{H}^{1/2}_{a,1}}+2\int_{t_1}^{t_2}\big(1-C_0\|u(s)\|_{\dot{H}^{1/2}_{a,1}}\big)\|u(s)\|^2_{\dot{H}^{3/2}_{a,1}}ds\le\|u(t_1)\|^2_{\dot{H}^{1/2}_{a,1}},
$$
so $t\mapsto\|u(t)\|_{\dot{H}^{1/2}_{a,1}}$ is nonincreasing on $[0,T_3)$ and $\|u(t)\|_{\dot{H}^{1/2}_{a,1}}\le\|u(0)\|_{\dot{H}^{1/2}_{a,1}}<\frac{1}{C_0}$ there. If $T_3<T^*$, continuity would give $\|u(T_3)\|_{\dot{H}^{1/2}_{a,1}}\le\|u(0)\|_{\dot{H}^{1/2}_{a,1}}<\frac{1}{C_0}$, and the strict inequality would persist after $T_3$, contradicting the definition of $T_3$. Hence $T_3=T^*$, and, using $1-C_0\|u(s)\|_{\dot{H}^{1/2}_{a,1}}\ge\kappa$, we get for all $t<T^*$
$$
\|u(t)\|^2_{\dot{H}^{1/2}_{a,1}}+2\kappa\int_0^{t}\|u(s)\|^2_{\dot{H}^{3/2}_{a,1}}ds\le\|u(0)\|^2_{\dot{H}^{1/2}_{a,1}} .
$$
In particular $\int_0^{T^*}\|u(t)\|^2_{\dot{H}^{3/2}_{a,1}}dt\le\frac{\|u(0)\|^2_{\dot{H}^{1/2}_{a,1}}}{2\kappa}<\infty$, and (ii) implies $T^*=+\infty$. $\qquad\Box$
 
\section{ Proof of Theorem~\ref{theo4}}
Let $u\in C(\mathbb{R}^+,\dot{H}^{1/2}_{a,1})$ be a global solution of (NS), that is, $u$ satisfies the integral equation $u=e^{t\Delta}u^0+B(u,u)$ on $[0,\infty)$, with $u^0=u(0)$. By the Remark following Theorem~\ref{theo1}, $u\in L^2_{loc}([0,\infty),\dot{H}^{3/2}_{a,1})$. We denote by $C_0$ the constant of Lemma~\ref{lem3} (it does not depend on $a$, so it is also valid for $a=0$), and by $\epsilon_0,\rho_0,c_0$ the universal constants of Lemma~\ref{theo3}. Fix
$$
0<\epsilon\le\min\Big\{\epsilon_0,\ \frac{1}{2C_0},\ \frac{\rho_0}{3(1+c_0)}\Big\}.
$$
We shall show that $\|u(t)\|_{\dot{H}^{1/2}_{a,1}}\le c_0\epsilon$ for all $t$ large enough; since $\epsilon$ is arbitrary, this proves (\ref{T4}).
 
\medskip\noindent{\bf Step 1: smallness at a time $t_0$.} Since $\dot{H}^{1/2}_{a,1}\hookrightarrow\dot{H}^{1/2}$, $u$ is a global solution in $C(\mathbb{R}^+,\dot{H}^{1/2})$, and by the result of Gallagher, Iftimie and Planchon \cite[Theorem~1]{GI} we have $\|u(t)\|_{\dot{H}^{1/2}}\to0$ as $t\to\infty$. Hence there is $t_0>0$ such that
$$
\epsilon':=\|u(t_0)\|_{\dot{H}^{1/2}}<\epsilon .
$$
Set $w(t):=u(t_0+t)$ for $t\ge0$. By the semigroup property, $w=e^{t\Delta}w(0)+B(w,w)$ on $[0,\infty)$, ${\rm div}\,w=0$, and $w(0)=u(t_0)$.
 
\medskip\noindent{\bf Step 2: $w$ belongs to the ball of radius $\rho_0$ of $X$.} We use the energy inequality (\ref{enqS2}) with $a=0$ (its proof only uses $w\in C(\dot{H}^{1/2})\cap L^2_{loc}(\dot{H}^{3/2})$). Since $C_0\epsilon'\le\frac12$, the continuity argument of the proof of (i) in Section~4 (with $a=0$ and $\kappa=1-C_0\epsilon'\ge\frac12$) applies and yields, for all $t\ge0$,
\begin{eqnarray}\label{enq51}
\|w(t)\|_{\dot{H}^{1/2}}\le\epsilon',\qquad \int_0^\infty\|w(s)\|^2_{\dot{H}^{3/2}}ds\le\frac{\epsilon'^2}{2\kappa}\le\epsilon'^2 .
\end{eqnarray}
Now $w$ is the solution given by Duhamel's formula of the heat equation with data $w(0)$ and source $f=-\mathbb{P}\,{\rm div}(w\otimes w)$. By (\ref{enq3}), the inequality $\|w\|^2_{\dot{H}^1}\le\|w\|_{\dot{H}^{1/2}}\|w\|_{\dot{H}^{3/2}}$ and (\ref{enq51}),
\begin{align*}
\|f\|^2_{L^2(\dot{H}^{-1/2})}&\le C_0^2\int_0^\infty\|w(s)\|^4_{\dot{H}^1}ds\\&\le C_0^2\,\epsilon'^2\int_0^\infty\|w(s)\|^2_{\dot{H}^{3/2}}ds\\&\le C_0^2\epsilon'^4 ,
\end{align*}
so $\|f\|_{L^2(\dot{H}^{-1/2})}\le C_0\epsilon'^2\le\frac{\epsilon'}{2}$. Lemma~\ref{lem4} (with $a=0$, on $[0,T]$, then $T\to\infty$) gives
$$
\|w\|_{\widetilde{L}^\infty(\dot{H}^{1/2})}\le\epsilon'+\frac{1}{\sqrt2}\cdot\frac{\epsilon'}{2}\le2\epsilon',
\qquad
\|w\|_{L^2(\dot{H}^{3/2})}\le\epsilon'
$$
(the second bound is (\ref{enq51})). Hence $w\in X=\widetilde{L}^\infty(\dot{H}^{1/2})\cap L^2(\dot{H}^{3/2})$ and
$$
\|w\|_X\le3\epsilon'<3\epsilon\le\frac{\rho_0}{1+c_0}<\rho_0 .
$$
 
\medskip\noindent{\bf Step 3: the analytic solution and uniqueness.} Since $\|w(0)\|_{\dot{H}^{1/2}}=\epsilon'<\epsilon\le\epsilon_0$ and ${\rm div}\,w(0)=0$, Lemma~\ref{theo3} gives a global solution $v\in X$ of
$$\begin{cases}
\partial_t v-\Delta v+v\cdot\nabla v=-\nabla p_1\\
{\rm div\,} v=0\\
v(0)=u(t_0),
\end{cases}$$
which is unique in the ball of $X$ of radius $\rho_0$ and satisfies
$$
\|e^{\sqrt{t}|D|}v\|_{\widetilde{L}^\infty(\dot{H}^{1/2})}\le c_0\,\epsilon'<c_0\epsilon .
$$
By Step 2, $w$ is also a solution in $X$ with $\|w\|_X<\rho_0$. By the uniqueness in Lemma~\ref{theo3}, $w=v$, that is, $v(t)=u(t_0+t)$ for almost every $t\ge0$, hence for every $t\ge0$ by continuity of $w$ in $\dot{H}^{1/2}_{a,1}$ (the bounds below are closed under limits in $t$ by Fatou's lemma in $\xi$).
 
\medskip\noindent{\bf Step 4: conclusion.} Let $t\ge a^2$. Then $\sqrt t\ge a$, so $e^{a|\xi|}\le e^{\sqrt t|\xi|}$ and
\begin{align*}
\|u(t_0+t)\|_{\dot{H}^{1/2}_{a,1}}&=\|e^{a|D|}v(t)\|_{\dot{H}^{1/2}}\\&=\|e^{(a-\sqrt t)|D|}e^{\sqrt t|D|}v(t)\|_{\dot{H}^{1/2}}\\&\le\|e^{\sqrt t|D|}v(t)\|_{\dot{H}^{1/2}}\\&\le c_0\,\epsilon .
\end{align*}
Therefore $\|u(t)\|_{\dot{H}^{1/2}_{a,1}}\le c_0\epsilon$ for all $t\ge t_0+a^2$. Since $\epsilon>0$ is arbitrary (up to the fixed upper bound above), $\lim_{t\to\infty}\|u(t)\|_{\dot{H}^{1/2}_{a,1}}=0$. $\qquad\Box$
 
\section{ Proof of Theorem~\ref{theo5}}
We follow the strategy of \cite[Theorem~3]{GI} (stability in $\dot{H}^{1/2}$; complete proofs are in \cite{GI2}) and \cite{JB} (Lei--Lin spaces), adapted to the spaces $\dot{H}^{1/2}_{a,1}$ with the constants of Section 2. Let $C_0$ be the constant of Lemma~\ref{lem3} and
$$
K:=216\,C_0^4,\qquad \Lambda:=\int_0^\infty\|u(z)\|^4_{\dot{H}^1_{a,1}}dz .
$$
 
\medskip\noindent{\bf Step 0: $\Lambda<\infty$.} By Theorem~\ref{theo4}, there is $t_0>0$ such that $\|u(t)\|_{\dot{H}^{1/2}_{a,1}}\le\frac{1}{2C_0}$ for all $t\ge t_0$. Then (\ref{enqS2}) gives, for $t\ge t_0$,
$\|u(t)\|^2_{\dot{H}^{1/2}_{a,1}}+\int_{t_0}^t\|u(s)\|^2_{\dot{H}^{3/2}_{a,1}}ds\le\|u(t_0)\|^2_{\dot{H}^{1/2}_{a,1}}$, so $\int_{t_0}^\infty\|u\|^2_{\dot{H}^{3/2}_{a,1}}<\infty$; on $[0,t_0]$ this integral is finite because $u\in L^2_{loc}(\dot{H}^{3/2}_{a,1})$. Since $M:=\sup_{t\ge0}\|u(t)\|_{\dot{H}^{1/2}_{a,1}}<\infty$ ($u$ is continuous and tends to $0$), (\ref{enqS0}) gives
$$
\Lambda\le M^2\int_0^\infty\|u(z)\|^2_{\dot{H}^{3/2}_{a,1}}dz<\infty .
$$
 
\medskip\noindent{\bf Step 1: the equation for the difference.} Let $v\in C([0,T^*),\dot{H}^{1/2}_{a,1})\cap L^2_{loc}([0,T^*),\dot{H}^{3/2}_{a,1})$ be the maximal solution of (NS) with $v(0)=v^0$, given by Theorem~\ref{theo1} (here ${\rm div}\,v^0=0$). We shall prove that $T^*=\infty$ when
\begin{eqnarray}\label{enq61}
\|v^0-u(0)\|^2_{\dot{H}^{1/2}_{a,1}}\le\frac{1}{32C_0^2}\,e^{-K\Lambda}.
\end{eqnarray}
Put $w=v-u$ and $w^0=v^0-u(0)$. Since $v\otimes v-u\otimes u=w\otimes w+u\otimes w+w\otimes u$, the function $w$ solves
$$
\partial_tw-\Delta w+\mathbb{P}\,{\rm div}\,F=0,\qquad F:=w\otimes w+u\otimes w+w\otimes u,\qquad w(0)=w^0,
$$
and $w\in C([0,T^*),\dot{H}^{1/2}_{a,1})\cap L^2_{loc}([0,T^*),\dot{H}^{3/2}_{a,1})$. As in the proof of (\ref{enqS1}) (Section 4), for $0\le t<T^*$,
$$
\|w(t)\|^2_{\dot{H}^{1/2}_{a,1}}+2\int_0^t\|w(s)\|^2_{\dot{H}^{3/2}_{a,1}}ds=\|w^0\|^2_{\dot{H}^{1/2}_{a,1}}-2\int_0^t\langle{\rm div}\,F,w\rangle_{\dot{H}^{1/2}_{a,1}}ds,
$$
and $|\langle{\rm div}\,F,w\rangle_{\dot{H}^{1/2}_{a,1}}|\le\|F\|_{\dot{H}^{1/2}_{a,1}}\|w\|_{\dot{H}^{3/2}_{a,1}}$ by the Cauchy--Schwarz inequality.
 
\medskip\noindent{\bf Step 2: differential inequality.} By Lemma \ref{lem2} with $s=t=1$, exactly as in the proof of Lemma~\ref{lem3}, $\|f\otimes g\|_{\dot{H}^{1/2}_{a,1}}\le C_0\|f\|_{\dot{H}^1_{a,1}}\|g\|_{\dot{H}^1_{a,1}}$. Hence, using (\ref{enqS0}) and (\ref{enqS0b}),
$$
\|w\otimes w\|_{\dot{H}^{1/2}_{a,1}}\le C_0\|w\|_{\dot{H}^{1/2}_{a,1}}\|w\|_{\dot{H}^{3/2}_{a,1}},\qquad
\|u\otimes w\|_{\dot{H}^{1/2}_{a,1}}+\|w\otimes u\|_{\dot{H}^{1/2}_{a,1}}\le2C_0\|u\|_{\dot{H}^1_{a,1}}\|w\|^{1/2}_{\dot{H}^{1/2}_{a,1}}\|w\|^{1/2}_{\dot{H}^{3/2}_{a,1}} .
$$
Therefore, for almost every $t$,
$$
\frac{d}{dt}\|w\|^2_{\dot{H}^{1/2}_{a,1}}+2\|w\|^2_{\dot{H}^{3/2}_{a,1}}\le2C_0\|w\|_{\dot{H}^{1/2}_{a,1}}\|w\|^2_{\dot{H}^{3/2}_{a,1}}+4C_0\|u\|_{\dot{H}^1_{a,1}}\|w\|^{1/2}_{\dot{H}^{1/2}_{a,1}}\|w\|^{3/2}_{\dot{H}^{3/2}_{a,1}} .
$$
By Young's inequality, $\alpha\beta\le\frac12\beta^{4/3}+\frac{27}{32}\alpha^4$ for $\alpha,\beta\ge0$; with $\alpha=4C_0\|u\|_{\dot{H}^1_{a,1}}\|w\|^{1/2}_{\dot{H}^{1/2}_{a,1}}$ and $\beta=\|w\|^{3/2}_{\dot{H}^{3/2}_{a,1}}$, the last term is at most
$\frac12\|w\|^2_{\dot{H}^{3/2}_{a,1}}+K\|u\|^4_{\dot{H}^1_{a,1}}\|w\|^2_{\dot{H}^{1/2}_{a,1}}$ (indeed $\frac{27}{32}\cdot256\,C_0^4=216\,C_0^4=K$). Consequently
\begin{eqnarray}\label{enq62}
\frac{d}{dt}\|w\|^2_{\dot{H}^{1/2}_{a,1}}+\Big(\frac32-2C_0\|w\|_{\dot{H}^{1/2}_{a,1}}\Big)\|w\|^2_{\dot{H}^{3/2}_{a,1}}\le K\|u\|^4_{\dot{H}^1_{a,1}}\|w\|^2_{\dot{H}^{1/2}_{a,1}} .
\end{eqnarray}
 
\medskip\noindent{\bf Step 3: continuity argument.} By (\ref{enq61}), $\|w^0\|_{\dot{H}^{1/2}_{a,1}}\le\frac{1}{4\sqrt2\,C_0}<\frac{1}{4C_0}$. Let
$$
T=\sup\Big\{t\in[0,T^*):\ \|w(s)\|_{\dot{H}^{1/2}_{a,1}}<\frac{1}{4C_0}\ \text{ for all } s\le t\Big\}>0 .
$$
For $t<T$ we have $2C_0\|w\|_{\dot{H}^{1/2}_{a,1}}\le\frac12$, so (\ref{enq62}) gives $\frac{d}{dt}\|w\|^2+\|w\|^2_{\dot{H}^{3/2}_{a,1}}\le K\|u\|^4_{\dot{H}^1_{a,1}}\|w\|^2$. Set $Z(t)=\|w(t)\|^2_{\dot{H}^{1/2}_{a,1}}+\int_0^t\|w\|^2_{\dot{H}^{3/2}_{a,1}}ds$; then $Z'\le K\|u\|^4_{\dot{H}^1_{a,1}}Z$, and Gronwall's lemma yields, for all $t\in[0,T)$,
\begin{eqnarray}\label{enq63}
\|w(t)\|^2_{\dot{H}^{1/2}_{a,1}}+\int_0^t\|w(s)\|^2_{\dot{H}^{3/2}_{a,1}}ds\le\|w^0\|^2_{\dot{H}^{1/2}_{a,1}}e^{K\Lambda}<\frac{1}{16C_0^2}.
\end{eqnarray}
If $T<T^*$, continuity would give $\|w(T)\|_{\dot{H}^{1/2}_{a,1}}\le\frac{1}{4\sqrt2\,C_0}<\frac{1}{4C_0}$, and the strict inequality would persist after $T$, contradicting the definition of $T$. Hence $T=T^*$ and (\ref{enq63}) holds on $[0,T^*)$.
 
\medskip\noindent{\bf Step 4: $T^*=\infty$.} By (\ref{enq63}) and Step 0,
$$
\int_0^{T^*}\|v\|^2_{\dot{H}^{3/2}_{a,1}}ds\le2\int_0^{T^*}\|w\|^2_{\dot{H}^{3/2}_{a,1}}ds+2\int_0^{\infty}\|u\|^2_{\dot{H}^{3/2}_{a,1}}ds<\infty,
$$
so Theorem~\ref{theo2}\,(ii) gives $T^*=+\infty$. Thus $v$ is global, and (\ref{enq63}) holds for all $t\ge0$, which is the stated estimate. $\qquad\Box$
 
\section{ Appendix}
In this section we prove the decay result stated in Remark~\ref{appendix}: if $u\in C(\mathbb{R}^+,\dot{H}^{1/2}_{a,1})$ is a global solution of (NS), then for every $s\ge\frac12$,
\begin{eqnarray}\label{enqA1}
\|u(t)\|_{\dot{H}^s}=o\big(t^{-\frac{s-\frac12}{2}}\big),\quad t\to+\infty .
\end{eqnarray}
 
\medskip\noindent{\bf Step 1: $\|u(t)\|_{\dot{H}^s}\to0$ for all $s\ge\frac12$.} Expanding the exponential, $e^{2a|\xi|}=\sum_{k\ge0}\frac{(2a)^k|\xi|^k}{k!}$, so
$$
\|u(t)\|^2_{\dot{H}^{1/2}_{a,1}}=\sum_{k=0}^\infty\frac{(2a)^k}{k!}\,\|u(t)\|^2_{\dot{H}^{\frac{1+k}{2}}},
$$
and every term is nonnegative. More directly, for $s\ge\frac12$ we have $|\xi|^{2s}\le C_{s,a}|\xi|e^{2a|\xi|}$ with $C_{s,a}=\sup_{x>0}x^{2s-1}e^{-2ax}<\infty$, hence
$$
\|u(t)\|_{\dot{H}^s}\le C_{s,a}^{1/2}\|u(t)\|_{\dot{H}^{1/2}_{a,1}}\underset{t\to\infty}{\longrightarrow}0
$$
by Theorem~\ref{theo4}. This gives (\ref{enqA1}) for $s=\frac12$ and the qualitative decay $o(1)$ for all $s\ge\frac12$. The rate in (\ref{enqA1}) comes from the analytic estimate of Lemma~\ref{theo3}, as follows.
 
\medskip\noindent{\bf Step 2: the rate.} Let $s>\frac12$, $\delta:=s-\frac12$, and $M_\delta:=\sup_{x>0}x^\delta e^{-x}=(\delta/e)^\delta$. Let $\epsilon>0$ be as in Section~5 (small enough), and let $t_0$, $v$ be given by Steps 1--3 of Section~5, so that $v(t)=u(t_0+t)$ for all $t\ge0$ and
\begin{align*}
\|e^{\sqrt t|D|}v(t)\|_{\dot{H}^{1/2}}&\le\|e^{\sqrt{t}|D|}v\|_{\widetilde{L}^\infty(\dot{H}^{1/2})}\\&\le c_0\,\epsilon .
\end{align*}
For $t>0$ and $x=\sqrt t|\xi|$ we have $|\xi|^\delta=t^{-\delta/2}x^\delta\le M_\delta\,t^{-\delta/2}e^{\sqrt t|\xi|}$. Therefore
\begin{align*}
\|v(t)\|_{\dot{H}^s}&=\Big(\int|\xi|^{2\delta}\,|\xi|\,|\widehat v(t,\xi)|^2d\xi\Big)^{1/2}\\&\le M_\delta\,t^{-\delta/2}\,\|e^{\sqrt t|D|}v(t)\|_{\dot{H}^{1/2}}\\&\le M_\delta\,c_0\,\epsilon\,t^{-\delta/2}.
\end{align*}
Let $T\ge2t_0$ and $t=T-t_0\ge T/2$. Then
$$
\|u(T)\|_{\dot{H}^s}=\|v(t)\|_{\dot{H}^s}\le2^{\delta/2}M_\delta\,c_0\,\epsilon\,T^{-\delta/2},
$$
that is, $\limsup_{T\to\infty}T^{\delta/2}\|u(T)\|_{\dot{H}^s}\le2^{\delta/2}M_\delta c_0\,\epsilon$. Since $\epsilon>0$ is arbitrary (the time $t_0$ depends on $\epsilon$), $\lim_{T\to\infty}T^{\delta/2}\|u(T)\|_{\dot{H}^s}=0$, which is (\ref{enqA1}). $\qquad\Box$
 
\medskip\noindent{\bf Acknowledgments.} I would like to thank Jamel Benameur for his kind advice and insightful comments.

\end{document}